\documentclass[10pt]{article}
\usepackage{amssymb}
\usepackage{amsfonts}
\usepackage{amsmath}
\usepackage{graphicx}
\usepackage{amsfonts}
\usepackage{amssymb}
\usepackage{epstopdf}
\usepackage{color}

\usepackage[section]{placeins}

\newcommand{\R}{{\mathbb R}}

\newtheorem{theorem}{Theorem}[section]

\newtheorem{proposition}[theorem]{Proposition}

\newtheorem{remark}[theorem]{Remark}

\begin{document}
\thispagestyle{empty}

\begin{center}

\huge { Stochastic epidemic SEIRS models with a constant latency period}

\vspace{.5cm}

\normalsize {\bf Xavier Bardina\footnote{X. Bardina is partially
supported by the grant MTM2015-67802-P from MINECO.}, Marco
Ferrante\footnote{M. Ferrante is partially
supported by the grant 60A01-8451 from Universit\`a di Padova.}
, Carles Rovira$^{*}$\footnote{C. Rovira is partially
supported by the grant MTM2015-65092-P from MINECO/FEDER, UE and by
Visiting Professor Program 2015 of the Universit\`a di Padova.}}

{\footnotesize \it $^1$ Departament de Matem\`atiques,
Universitat Aut\`onoma de Barcelona,
08193 Bellaterra.

$^2$ Dipartimento di Matematica,
Universit\`a
di Padova,
Via Trieste 63, 35121-Padova, Italy.

$^3$ Departament de Matem\`atiques i Inform\`atica, Universitat de Barcelona, Gran
Via 585, 08007-Barcelona.

 {\it E-mail addresses}: Xavier.Bardina@uab.cat, ferrante@math.unipd.it,
Carles.Rovira@ub.edu}

{$^{*}$corresponding author}

\end{center}

\begin{abstract}
In this paper we consider the stability of a class of deterministic
and stochastic SEIRS epidemic models with delay. Indeed, we assume
that the transmission rate could be stochastic and the presence of a
latency period of $r$ consecutive days, where $r$ is a fixed
positive integer, in the ``exposed'' individuals class E. Studying
the eigenvalues of the linearized system, we obtain conditions for
the stability of the free disease equilibrium, in both the cases of
the deterministic model with and without delay. In this latter case,
we also get conditions for the stability of the coexistence
equilibrium. In the stochastic case we are able to derive a
concentration result for the random fluctuations and then, using the
Lyapunov method, that under suitable assumptions the free disease
equilibrium is still stable.
\end{abstract}



{\bf AMS 2000 MSC:} 92D30, 60J10, 60H10.

{\bf Running head:} A SEIRS stochastic model

\section{Introduction}
The mathematical models developed to describe the spread of a communicable disease
are both deterministic and stochastic and
they may involve many factors such as infectious agents, mode of transmission,
incubation periods, infectious periods, susceptibility, etc...

A well known deterministic model in a closed population consisting
of susceptible (S), infective (I) and recover (R) were considered
by Kermack and McKendrick in \cite{KM}.
Since then various epidemic deterministic models have been developed, such as SIR, SIS, SEIR and
SEIRS models with or without a time delay (see e.g McCluskey \cite{Mc1}, \cite{Mc2} and
Huang et al. \cite{HTMW}).
Here, the class E denotes individuals ``exposed'' to the disease, but not yet infectious.

\bigskip

In this paper we assume that our epidemic model encompasses the
class E. This is the case of several diseases like chickenpox
(discrete-time version of this model has been considered in
\cite{FFR}). Contrary to most of the models defined in the
literature, we consider a constant holding time in the E class.
i.e., we will assume a latency period of $r$ consecutive days, where
$r$ is a fixed positive integer. On the contrary, the permanence in
the other classes is defined in a classical manner. Finally, we
assume that recover individuals may became susceptible again and for
this reason our model will be a SEIRS model with delay.

More precisely, we will assume that the total size of the population will be fixed and equal to $n$ and that
all the individuals belong to one of the four classes S, E,  I or R, where
S denotes a susceptible individual, E an infected but not infectious individual (latency period),
I an infected and infectious individual and R a recovered individual.
These assumptions lead to the following set of ordinary differential equations, where
$S(t)$ denotes the fraction of individuals that are susceptible and the same for $E(t)$,
$I(t)$ and $R(t)$:
\begin{equation}
\label{SEIRS_d}
\left\{
\begin{array}{ll}
dS(t) = & -\beta S(t) I(t) dt + \gamma R(t) dt
\\
dE(t) = & \beta S(t) I(t) dt - \frac{1}{K_r} E(t-r) dt
\\
dI(t) = &\frac{1}{K_r} E(t-r) dt - \mu I(t) dt
\\
dR(t) = & \mu I(t) dt - \gamma R(t) dt
\end{array}
\right.
\end{equation}
where $\beta$ represent the disease transmission coefficient,
$\mu$ the rate at which infectious individuals becomes recovered,
$\gamma$ the rate at which recovered individuals become
again susceptible and $K_r$ the rate of latency. Note that $K_r$ is a parameter that depends on the delay $r$ on $E$.
The nonstandard equation in the model (\ref{SEIRS_d}) is represented by the second
one, since we assume that the change in the percentage of  individuals in the class E
depends on the difference between the number of individuals that enter
into this class at time $t$ and those that entered $r$ units of time before.

Similar deterministic models have been considered in the literature. Huang {\it et al}. \cite{HTMW}
consider the stability of  SIR and SIRS models with constant time delay caused by latency in a host
and expressed as a function of $S(t-r)$ and $I(t-r)$.  In Huang {\it et al}.  \cite{HBT} the authors consider
a SEIR model with constant latency  time and infectious periods. Bai \cite{B}  considers a delayed SEIRS
model with varying total population size where the delay  is also expressed as a function of $S(t-r)$ and $I(t-r)$.
Constant delays in epidemic models also appears in \cite{DZ}  to model relapse in infectious diseases.

\bigskip

We also consider a stochastic behaviour of the disease transmission, assuming that
$\beta$ may be random. To do this, we will assume that a infectious individual
makes a random number
$$\beta dt + \varepsilon d W_t$$
of contacts with other individuals in a time interval $[t,t+dt)$,
where $\{W_t, t \ge 0 \}$ denotes a standard Brownian motion
(see \cite{GGHMP} p.879).
Thus, the deterministic SEIRS model (\ref{SEIRS_d})  becames stochastic in the following way
\begin{equation}
\label{SEIRS_s}
\left\{
\begin{array}{ll}
dS(t) = & [-\beta S(t) I(t) + \gamma R(t)] dt -\varepsilon S(t) I(t) dW_t
\\
dE(t) = & [\beta S(t) I(t) - \frac{1}{K_r} E(t-r) ]dt + \varepsilon S(t) I(t) dW_t
\\
dI(t)= &[\frac{1}{K_r} E(t-r) - \mu I(t) ]dt
\\
dR(t) = &[ \mu I(t) - \gamma R(t)]dt
\end{array}
\right.
\end{equation}

As far as we know, the literature of epidemic stochastic models is scarce. Nevertheless, we may cite
Tornatore {\it et al}  \cite{Tor}
who proposes a stochastic SIR model with distributed time delay and discuss its stability and
Gray {\it et al.} \cite{GGHMP} who presents the study of a SIS epidemic model.

\bigskip

Let us point that the presence of a constant delay leads to (stochastic) delay differential equations,
which are not easy to handle mathematically.
Usually in the literature, the constant time delay appears with some incidence rate, i.e.,
$F(S(t-r))G(I(t-r))$ or with a control factor $e^{-kr}S(t-r)I(t-r)$.
These factors help to deal with the delay.
On the contrary, in our paper, we deal with a fixed delay in the class
$E$, that is, $E(t-r)$, since we consider that the number of individuals
that pass from $E$ to $I$ al time $t$, depends only on the number of
individuals at $E$ at time $t-r$.
We also have a coefficient $K_r^{-1}$ that depends on $r$ but the relation
is not exponentially and it has been considered in order to discuss the validity and
stability of the model depending on the value of $r$.

\bigskip

In Section \ref{det} we study the deterministic model
(\ref{SEIRS_d}). We analyse first the case without delay $r=0$
obtaining that when $\mu \ge \beta$  the free disease equilibrium is
stable while when  $\beta > \mu$ the stability holds in the point of
coexistence equilibrium. Thus we can see that it does not depend on
$\gamma$ (notice that $\gamma$ allows that some recovered
individuals became again susceptibles).  Then we study the delayed
model. We get that to ensure the validity of the model it is
necessary that $K_r \ge re$. Then, we show that  when $\mu \ge
\beta$  the free disease equilibrium remains assymptotically stable
for any delay $r$.

Section \ref{stoc} is devoted to deal with the stochastic model
(\ref{SEIRS_s}). Our aim is to study what happens with the stochastic fluctuations of the deterministic model. We show that the solutions of the perturbed system tends uniformly (and exponentially) to the solutions of the deterministic model when
$\varepsilon$ tends to zero. It holds for both models, with delay and without delay. Finally, we show that for the nondelayed model, under the condition
\begin{equation*}
\mu - \beta - \frac{1}{\mu K_r}\frac{\varepsilon^2}{2} >0.
\end{equation*}
the free disease equilibrium is assymptotically stable.

Finally, in the Appendix we recall some basic results about stability: methods based on the study of the roots of the associated characteristic functions (we use them in Section \ref{det}) and methods based on Lyauponov functionals (used in Section \ref{stoc}).

\section{Deterministic model}\label{det}

Let us now consider the model
(\ref{SEIRS_d})
with the initial condition $E(s)=e_0, S(t)+I(t)+R(t)=0$ for $s\in(-r,0]$. In the first subsection we will  deal with the model without delay, i.e. when $r=0$. In the second subsection we will study what happens when we introduce the delay.

\subsection{Analysis of the case without delay}

\subsubsection{Existence and Positivity of the Solution}

Using standard method we get the existence and uniqueness of
solution. It follows from the fact that if we start with $s_0 \ge 0,
e_0 \ge 0, i_0 \ge 0 $ and $r_0 \ge 0$,  the region
$$
\{ (s,e,i,r): s,e,i,r \ge 0; s+e+i+r \le 1 \}$$ is positively
invariant.

 For instance, if $S(t_1)=0, E(t_1)>0, I(t_1)>0,R(t_1)>0$  for some $t \ge 0$ then  $\dot S(t_1)=\gamma R(t_1) >0$ and there will exists $\epsilon$ such that $S(t)>0$  for any $t \in (t_1, t_1 + \epsilon)$.  All the other cases can be done by similar arguments. On the other hand, obviously $S(t)+E(t)+I(t)+R(t) \le 1$.

\subsubsection{Analysis of the equilibrium}

To get equilibria, we have to compute the solution to the following equations:
\begin{equation}
\label{equil}
\left\{
\begin{array}{ll}
0 = & -\beta S  I  + \gamma R
\\
0 = & \beta S  I  - \frac{1}{K_r} E
\\
0= & \frac{1}{K_r}  E  - \mu I \\
0 = & \mu I  - \gamma R
\end{array}
\right.
\end{equation}
with the restriction $S+E+I+R=1$. We have that the free disease equilibrium $X^0=(1,0,0,0)$ exists for any
values of the parameters. In the case $\mu \ge \beta$ one can check that no other equilibrium exists while in the case
$\beta > \mu$ there exists also one point of coexistence equilibrium
 {\small
\begin{equation*}X^{*}=\left({\frac {\mu}{\beta}},{\frac {K_r \left( \beta-\mu
\right) \mu\,\gamma}{ \beta\, \left( \gamma\,K_r\mu+\gamma+\mu
\right) }},{\frac {\gamma\,
 \left( \beta-\mu \right) }{\beta\, \left( \gamma\,K_r\mu+\gamma+\mu
 \right) }},{\frac { \left( \beta-\mu \right) \mu}{\beta\, \left(
\gamma\,K_r\mu+\gamma+\mu \right) }}\right)\end{equation*} }

\begin{proposition}\label{propa}

\begin{enumerate}

\item If $\mu \ge \beta$ , then $X^0$ is  asymptotically stable;

\item If $\mu < \beta$, then $X^*$ exists and is  asymptotically stable.

\end{enumerate}
\end{proposition}

\begin{remark}
It can be checked easily that the basic reproduction number for this model is
$$R_0=\frac{\beta}{\mu}.$$ So, the stability of the free disease equilibrium hold when $R_0 \le 1$ while otherwise we have the stability
of the coexistence equilibrium.
\end{remark}

\noindent{\bf Proof of Proposition \ref{propa}:}
Since we have the relation $S(t)=1-E(t)-I(t)-R(t)$, we can consider that we are dealing with the 3-dimensional system
\begin{equation}
\label{SEIRS_2}
\left\{
\begin{array}{ll}
dE(t) = & [\beta (1-E(t)-I(t)-R(t)) I(t) - \frac{1}{K_r} E(t-r)]dt
\\
dI(t)= &[\frac{1}{K_r} E(t-r) - \mu I(t)]dt
\\
dR(t)= & [\mu I(t) - \gamma R(t)]dt
\end{array}
\right.
\end{equation}
The coefficient matrix of the linearized system at the free disease equilibrium is
{\[\displaystyle  \left( \begin {array}{ccc}
\noalign{\medskip}-{{\it K_r}}^{-1}&\beta&0\\
\noalign{\medskip}{{\it K_r}}^{-1}&-\mu&0\\
\noalign{\medskip}0&\mu&-\gamma\end {array} \right)  \]}
with eigenvalues
 {\[\displaystyle  \left[ \begin
{array}{c} \noalign{\medskip}\frac12\,{\frac {-\mu\,{\it
K_r}-1+\sqrt {{\mu}^{2}{{\it K_r}}^{2}-2\,\mu\,{\it K_r}+4\,{\it
K_r}\,\beta+1}}{{\it K_r}}}\\ \noalign{\medskip}-\frac12\,{\frac
{\mu\,{\it K_r}+1+\sqrt {{\mu}^{2}{{\it K_r}}^{2}-2\,\mu\,{\it
K_r}+4\,{\it K_r}\,\beta+1}}{{\it K_r}}}\\
\noalign{\medskip}-\gamma\end {array} \right] .\]}
The last two are clearly negative. On the other hand, we have
\begin{eqnarray*}&&\frac12\,{\frac {-\mu\,{\it
K_r}-1+\sqrt {{\mu}^{2}{{\it K_r}}^{2}-2\,\mu\,{\it K_r}+4\,{\it
K_r}\,\beta+1}}{{\it K_r}}}\\
&=& \frac12\,{\frac {-(\mu\,{\it K_r}+1)+\sqrt {({\mu}{{\it
K_r}}+1)^{2}-4\,{\it K_r}(\mu-\beta)}}{{\it K_r}}}
\end{eqnarray*}
that is also negative when  $\beta<\mu$. Thus, when   $\beta<\mu$ all the eigenvalues are negative
and so, the free disease equilibrium is locally asymptotically estable.

Let us consider now what happens around the coexistence equilibrium.
The coefficient matrix of the linearized system at the coexistence
equilibrium is now {\[\displaystyle A:=\left( \begin {array}{ccc}
-{\frac {\beta\,\gamma\,K_r+\gamma+\mu}{K_r
 \left( \gamma\,K_r\mu+\gamma+\mu \right) }}&{\frac {\gamma\,K_r{\mu}^{2}-
\beta\,\gamma+2\,\gamma\,\mu+{\mu}^{2}}{\gamma\,K_r\mu+\gamma+\mu}}&-{
\frac {\gamma\, \left( \beta-\mu \right)
}{\gamma\,K_r\mu+\gamma+\mu}}
\\ \noalign{\medskip}{K_r}^{-1}&-\mu&0\\ \noalign{\medskip}0&\mu&-\gamma
\end {array} \right).
\]}
Our aim is to chek that when $\beta > \mu$ all the eigenvalues have
negative real part. Using the well-known Routh-Hurwitz criterium, it
is enough to chek that $\textrm{Trace}(A)<0$,
$\textrm{Determinant}(A)<0$ and
$-A_2*\textrm{Trace}(A)+\textrm{Determinant}(A)>0$ where $A_2$ is
the coefficient of $\lambda$ in the characteristic polynomial
$P(\lambda)$ of $A$, i.e., if $A=(a_{i,j})$,
$$A_2:=a_{1,1}a_{2,2}+a_{1,1}a_{3,3}+a_{2,2}a_{3,3}-a_{1,2}a_{2,1}-a_{1,3}a_{3,1}-a_{2,3}a_{3,2}.$$
Indeed
\begin{eqnarray*}
\textrm{Trace($A$)}=-{\frac {\beta\,\gamma\,K_r+\gamma+\mu}{K_r
\left( \gamma\,K_r\mu+\gamma+\mu
 \right) }}-\mu-\gamma<0\\
\textrm{Determinant($A$)}=-{\frac {\gamma\, \left( \beta-\mu \right)
}{K_r}}<0,\end{eqnarray*} and finally
$$A_2={\frac {\gamma\, \left(
\gamma\,{K_r}^{2}{\mu}^{2}+\beta\,\gamma\,K_r+
\beta\,K_r\mu+\gamma\,K_r\mu+\beta+\gamma \right) }{K_r \left(
\gamma\,K_r\mu+ \gamma+\mu \right) }}.$$
Notice that,
$$
A_2>{\frac {\gamma\, \left( \beta\,K_r\mu+\beta\right) }{K_r \left(
\gamma\,K_r\mu+ \gamma+\mu \right) }}, \quad
-\textrm{Trace($A$)}>\mu+\gamma,
$$
and
$$
\textrm{Determinant($A$)}>-{\frac {\gamma\, \beta}{K_r}}.
$$
So
\begin{eqnarray*}
&&-A_2*\textrm{Trace($A$)}+\textrm{Determinant($A$)}\\&>&{\frac
{\gamma\, \left( \beta\,K_r\mu+\beta\right) }{K_r \left(
\gamma\,K_r\mu+ \gamma+\mu \right) }}(\mu+\gamma)-{\frac {\gamma\,
\beta}{K_r}}=\frac{\gamma\,\beta\,K_r\,\mu^2}{K_r \left(
\gamma\,K_r\mu+ \gamma+\mu \right) }>0.\end{eqnarray*}

\hfill$\square$

\subsection{Analysis of the delayed case}

\subsubsection{Existence and Positivity of the solution}

The system (\ref{SEIRS_d}) can be solved step by step.
Indeed, if we are able to solve the system up to time $n r$, we can find
the solution for $t\in [n r,(n+1)r]$, since
\[
I(t)= e^{-\mu(t-n r)}
\left(\frac{1}{K_r}
\int_{n r}^{t}  E(s-r)e^{\mu(s-n r)} ds + I(n r)
\right)
\]
and moreover
\begin{eqnarray*}
R(t)&= &e^{-\gamma(t-n r)}
\left(
\int_{n r}^{t} \mu  I(s)e^{\gamma(s-n r)} ds + R(n r) \right) \\
S(t)
&=& e^{-\int_{n r}^{t} \beta I(s) ds}
\left(
\int_{n r}^{t} \gamma  R(s)e^{\int_{n r}^{s} \beta I(u) du}  ds + S(n r) \right) \\
E(t)&= & E(nr) + \int_{n r}^{t} \beta S(s)I(s) ds -\frac{1}{K_r} \int_{n r}^{t}  E(s-r) ds.
\end{eqnarray*}
On the other hand, we can reduce the problem to the study of the positivity of $E$, since if $E$ is nonnegative on $[0,t]$ then $I, R$ and $S$ are clearly nonnegative functions
 on $[0,t+r]$.  Moreover, if $I$ and $S$ are nonnegative on  $[0,t+r]$ we have that
$$
\frac{dE}{ds} = \beta S(s) I(s) - \frac{1}{K_r} E(s-r) \ge  -  \frac{1}{K_r} E(s-r) $$
and using a comparision argument $E(s) \ge F(s)$ where $F$ is the solution of the equation
$$
\frac{dF}{ds} =  -  \frac{1}{K_r} F(s-r) .$$

The next remark helps us to justify the positivity of the solutions.
\begin{remark}
Let us consider the delay differential equation
$$
\frac{dF}{ds} =  -  K F(s-r) $$
on $[0,t]$ where $K$  is apositive constant.  From Theorem A in \cite{PK} it follows that if $K \le \frac{1}{r} e^{-1}$ then the solution of this equation is nonnegative. Furthermore,  Corollary 2.1 in \cite{Dor} yields that if $K = \frac{1}{r} e^{-1}$ then the solution is nonnegative and $\lim_{s \to \infty} F(s)=0$.
\end{remark}

So, we can state the existence and positivity in the following proposition.

\begin{proposition}
Assume that  $\beta, \mu, \gamma \in (0,1)$  and  $K_r \ge re$, then the system (\ref{SEIRS_d}) has an unique nonnegative solution.
\end{proposition}


\subsubsection{Analysis of the equilibrium points}

The linearization of the delayed SEIR system about the steady state (1,0,0,0) is

{\small {\[\displaystyle
 \left(\!\!\! \begin {array}{c}
\noalign{\medskip}\frac{dS}{dt}\\
\noalign{\medskip}\frac{dE}{dt}\\
\noalign{\medskip}\frac{dI}{dt}\\
\noalign{\medskip}\frac{dR}{dt}\end {array} \!\!\!\right) =
 \left( \begin {array}{cccc}
\noalign{\medskip}0&0&-\beta&+\gamma\\
\noalign{\medskip}0&0&\beta&0\\
\noalign{\medskip}0&0&-\mu&0\\
\noalign{\medskip}0& 0&\mu&-\gamma\end {array} \right)
 \left(\!\!\! \begin {array}{c}
\noalign{\medskip}S(t)\\
\noalign{\medskip}E(t)\\
\noalign{\medskip}I(t)\\
\noalign{\medskip}R(t)\end {array} \!\!\!\right) +
 \left( \begin {array}{cccc}
\noalign{\medskip}0&0&0&0\\
\noalign{\medskip}0&-{{\it K_r}}^{-1}&0&0\\
\noalign{\medskip}0&{{\it K_r}}^{-1}&0&0\\
\noalign{\medskip}0& 0&0&0\end {array} \right)
 \left(\!\!\! \begin {array}{c}
\noalign{\medskip}S(t-r)\\
\noalign{\medskip}E(t-r)\\
\noalign{\medskip}I(t-r)\\
\noalign{\medskip}R(t-r)\end {array} \!\!\!\right).
\]}}
As in the nondelayed case we can eliminate  the first equation since $S(t)=1-(E(t)+I(t)+R(t))$. Thus we get the
 characteristic equation
$$
P(\lambda)=(\lambda + \gamma) \Big( \lambda^2 + \lambda \mu + \big( \lambda \frac{1}{K_r} + (\frac{\mu}{K_r}- \frac{\beta}{K_r})\big)e^{-\lambda r} \Big).
$$
Let us fix $K_r$ and let us study what happens when the delay $r$ is increasing if $\beta
< \mu$, .  We have to study the behaviour of the roots of the characteristic function. Clearly $-\gamma $ is a negative real root. So, our aim will be the study of the other factor
$$\lambda^2 + \lambda \mu + \big( \lambda \frac{1}{K_r} + (\frac{\mu}{K_r}- \frac{\beta}{K_r})\big)e^{-\lambda r}=0.$$
Applying Proposition \ref{Prop2}  with $a_1=\mu, a_0=0, b_1=K_r^{-1}, b_0=
K_r^{-1}(\mu - \beta)$, it is easy to check that $a_0 + b_0 >0, a_1+b_1>0,$ and
 $a_0^2 < b_0^2$ and consequently there will be positive roots for some $r$ and  the steady sate $(1,0,0,0)$ became unstable when $r$ is increasing. 

Following the results in Subsection \ref{secdeg2} and applying (\ref{arrel1}), (\ref{eqcos}) and (\ref{eqsin}),  the state (1,0,0,0) remains asymptotically stable until $r^*=\frac{\theta}{\omega}$, where $\omega>0,  0 \le \theta <2 \pi$ and
\begin{eqnarray*}
\omega^2 &= &\frac12  \Big( (K_r^{-2}-\mu^2) + \big((K_r^{-2}-\mu^2)^2+4 K_r^{-2}(\mu - \beta)^2\big)^\frac12 \Big),\\
\cos \theta & = & - \frac{\mu K_r^{-1} \omega^2  - \omega^2 (\mu-\beta)  K_r^{-1} }{ K_r^{-1} \omega^2 + K_r^{-2}(\mu - \beta)^2}=
 - \frac{\beta K_r^{-1} \omega^2  }{ K_r^{-1} \omega^2 + K_r^{-2}(\mu - \beta)^2} , \\
\sin \theta & = & \frac{ K_r^{-1} (\mu - \beta) \mu \omega +K_r^{-1} \omega^3 }{ K_r^{-1} \omega^2 + K_r^{-2}(\mu - \beta)^2}.
\end{eqnarray*}
Since $\beta < \mu$, we clearly have that $\cos \theta <0$ and $ \sin \theta >0$ and furthermore that $\theta \ge \frac{\pi}{2}.$
Thus, if
$$r \le  M(K_r,\mu,\beta) := \frac{\pi}{2^\frac12 \Big(  (K_r^{-2}-\mu^2) + \big((K_r^{-2}-\mu^2)^2+4 K_r^{-2}(\mu - \beta)^2\big)^\frac12  \Big)^\frac12}.$$
  the state (1,0,0,0) remains assymptotically stable.
Notice that
\begin{eqnarray*}
M(K_r,\mu,\beta) & \ge  &M(K_r,\mu,0) \\&=&
\frac{\pi}{2^\frac12 \Big(  (K_r^{-2}-\mu^2) + \big((K_r^{-2}-\mu^2)^2+4 K_r^{-2}\mu^2\big)^\frac12  \Big)^\frac12}\\
    &=&\frac{1}{2} \pi K_r.
\end{eqnarray*}
Thus, since  $r < \frac{K_r}{e},$ it holds that $r < \frac{1}{2} \pi K_r \le M(K_r,\mu,\beta)$ and
the state $(1,0,0,0)$ remains asympotically stable for any possible
delay. The result states as follows:

\begin{proposition}
Assume that  $\beta, \mu, \gamma \in (0,1)$  and  $K_r \ge re$, then the free disease equilibium point is asymptotically stable for any $r>0$.
\end{proposition}


\bigskip

Let us study now what happens with the coexistence equilibrium when $\beta > \mu$.  The linearization of the delayed SEIR
system about the coexistence equilibrium, reduced to three
equations, is
\begin{eqnarray*}
 \left( \begin {array}{c}
\noalign{\medskip}\frac{dE}{dt}\\
\noalign{\medskip}\frac{dI}{dt}\\
\noalign{\medskip}\frac{dR}{dt}\end {array} \right) &=&
 \left( \begin {array}{ccc} -{\frac { \left( \beta-\mu \right) \gamma
}{\gamma\,K_r\mu+\gamma+\mu}}&-{\frac { \left( \beta-\mu \right)
\gamma }{\gamma\,K_r\mu+\gamma+\mu}}+\mu&-{\frac { \left( \beta-\mu
\right) \gamma}{\gamma\,K_r\mu+\gamma+\mu}}\\
\noalign{\medskip}0&-\mu&0
\\ \noalign{\medskip}0&\mu&-\gamma\end {array}
 \right)
 \left( \begin {array}{c}
\noalign{\medskip}E(t)\\
\noalign{\medskip}I(t)\\
\noalign{\medskip}R(t)\end {array} \right) \\&&+
 \left( \begin {array}{ccc}
\noalign{\medskip}-{{\it K_r}}^{-1}&0&0\\
\noalign{\medskip}{{\it K_r}}^{-1}&0&0\\
\noalign{\medskip}0&0&0\end {array} \right)
 \left( \begin {array}{c}
\noalign{\medskip}E(t-r)\\
\noalign{\medskip}I(t-r)\\
\noalign{\medskip}R(t-r)\end {array} \right).
\end{eqnarray*}
As in the previous case we get the
 characteristic equation
\begin{equation}\label{deg33}
\lambda^3 + a_2 \lambda^2 + a_1 \lambda + a_0 + \big( b_2 \lambda^2 +b_1 \lambda  +b_0 \big)e^{-\lambda r}=0.
\end{equation}
with
\begin{eqnarray*}
{\it a_0}&=&{\frac {{\gamma}^{2}\mu\, \left( \beta-\mu \right)
}{\gamma \,K_r\mu+\gamma+\mu}},
\\
{\it a_1}&=&{\frac {\gamma\, \left(
\gamma\,K_r{\mu}^{2}+\beta\,\gamma+ \beta\,\mu \right)
}{\gamma\,K_r\mu+\gamma+\mu}},
\\
{\it a_2}&=&{\frac
{{\gamma}^{2}K_r\mu+\gamma\,K_r{\mu}^{2}+\beta\,\gamma+{
\gamma}^{2}+\gamma\,\mu+{\mu}^{2}}{\gamma\,K_r\mu+\gamma+\mu}},
\\
{\it b_0}&=&{\frac {\gamma\, \left(
\beta\,\gamma+\beta\,\mu-\gamma\, \mu-{\mu}^{2} \right) }{ \left(
\gamma\,K_r\mu+\gamma+\mu \right) K_r}},
\\{\it b_1}&=&{\frac {\gamma\, \left( \gamma\,K_r\mu+\beta+\gamma \right)
}{ \left( \gamma\,K_r\mu+\gamma+\mu \right) K_r}},
\\ {\it b_2}&=&{{\it
K_r}}^{-1}.
\end{eqnarray*}
Applying Proposition \ref{Prop24} to our characteristic function we
get that the state will become unstable with increasing delay since it is easy to check that
 $a_2+b_2>0, a_0+b_0>0,(a_2+b_2)(a_1+b_1)-(a_0+b_0)>0$ and, with the additional condition $K_r<\frac1\mu+\frac1\gamma$, $C=a_0^2 - b_0^2<0$.

Although this case is more complicated and   it is not possible to give a general
criterium as in the endemic case, following the methods in Kuang in \cite{Ku}
pages 74-76, we  will find some conditions about  the first crossing of the imaginary axis.
That is, we will find a point to which the steady state will remain asymptotically stable.

Let us assume that $\lambda = i \omega, \omega>0$, is a root for some $r$.
Under the  hypothesis $a_0+b_0 \not= 0$, we clearly have that $\omega \not= 0$.  Then, we have
\begin{eqnarray*}
(a_0-a_2 \omega^2) + b_1 \omega \sin(\omega r) + (b_0-b_2 \omega^2) \cos(\omega r) & = & 0 \\
(- \omega^3+a_1 \omega) + b_1 \omega \cos(\omega r) + (b_0-b_2 \omega^2) \sin(\omega r) & = & 0.
\end{eqnarray*}
So
\begin{eqnarray}
0 & = & (a_2^2 \omega^2 - a_0)^2 + (\omega^3 - a_1 \omega)^2 - b_1^2 \omega^2 - (b_0-b_0\omega^2)^2 \nonumber\\
& =& \omega^6 + A \omega^4 + B \omega^2 +  C, \label{garu3}
\end{eqnarray}
where $A, B$ and $C$ are defined in (\ref{ABC}).
Let us consider now the associated discriminant to the third degree equation (\ref{garu3})
\begin{equation*}
\Delta:= 18ABC-4A^3C+A^2B^2 -4B^2 -27 C^2.
\end{equation*}
Then if $\Delta <0$, the equation (\ref{garu3})  (as an equation of $\omega^2$) has only one real root $\omega_0$ and the system will remain stable unless until the delay
$r=\frac{\theta}{\omega_0}$ where
\begin{eqnarray} \label{Marco}
\cos \theta & = & - \frac{ (- \omega_0^3+a_1  \omega_0)b_1 \omega_0 + (a_0 - a_2\omega_0^2) (b_0-b_2 \omega_0^2) }{ b_1^2 \omega_0^2 +( b_0-b_2 \omega_0^2)^2}  \\
\sin \theta & = & \frac{ (b_0-b_2 \omega_0^2)  (- \omega_0^3+a_1  \omega_0)  - b_1 \omega_0  (a_0 - a_2\omega_0^2) }{ b_1^2 \omega_0^2 +( b_0-b_2 \omega_0^2)^2}. \nonumber
\end{eqnarray}

This leads to the following result:

\begin{proposition}
Assume that  $\beta, \mu, \gamma \in (0,1)$  and  $K_r \ge re$ with $K_r<\frac1\mu+\frac1\gamma$, then if $\Delta<0$ the  coexistence  equilibium point is asymptotically stable for any $r<\frac{\theta}{\omega_0}$ where $\theta$ and $\omega_0$ satisfy (\ref{Marco}).
\end{proposition}

\section{Stochastic model }\label{stoc}

In this section we get some stability properties for the stochastic system.
More precisely we will study the fluctuations of the random system. Let us recall that we are considering the system:
\begin{equation*}
\left\{
\begin{array}{ll}
dS^\varepsilon(t) = & (-\beta S^\varepsilon(t) I^\varepsilon(t) + \gamma R^\varepsilon(t)) dt + \varepsilon S^\varepsilon(t) I^\varepsilon(t) dW(t)
\\
dE^\varepsilon(t)= & (\beta S^\varepsilon(t) I^\varepsilon(t) - \frac{1}{K_r} E^\varepsilon(t-r) )dt - \varepsilon S^\varepsilon(t) I^\varepsilon(t) dW(t)
\\
dI^\varepsilon(t)= &(\frac{1}{K_r} E^\varepsilon(t-r) - \mu I^\varepsilon(t) ) dt
\\
dR^\varepsilon(t) = & (\mu I^\varepsilon(t) - \gamma R^\varepsilon(t)) dt
\end{array}
\right.
\end{equation*}

We get a concentration result for the random fluctuations, that holds for the delayed and the non-delayed systems. This exponential stability states as follows.

\begin{proposition}
Assume that  $\beta, \mu, \gamma \in (0,1)$  and  $K_r \ge re$.  Set $Z^\varepsilon(t) $ for the random vector $(S^\varepsilon(t),E^\varepsilon(t),I^\varepsilon(t),R^\varepsilon(t))$ and $Z(t)$ for the solution to the corresponding deterministic system.  Then, there exists  nonnegative constants $K_1, K_2$ depending on $r,\beta, \mu, \gamma$ such that
$$
P( \Vert Z^\varepsilon-Z \Vert_{\infty,[0,T]} > \rho )  \le
\exp \Big( - \frac{\rho^2}{\varepsilon^2 K_1 T \exp(K_2 T)} \Big) .
$$
\end{proposition}

\noindent{\bf Proof:}
Set $J^\varepsilon(t):=\varepsilon \int_0^t S^\varepsilon(s) I^\varepsilon(s) dW(s).$ We can write, using that $S^\varepsilon$ and $I$ are bounded:

\begin{eqnarray*}
| S^\varepsilon(t)-S(t) | & \le & \beta \int_0^t | S^\varepsilon(u) I^\varepsilon(u)-S(u) I(u)| du \\ & &  \qquad+ \int_0^t \gamma | R^\varepsilon(u) -R(u)| du
+| J^\varepsilon(t) |\\
& \le & \beta \int_0^t | S^\varepsilon(u)-S(u)| du  + \beta \int_0^t
|  I^\varepsilon(u)- I(u)| du\\ & &  \qquad+ \gamma \int_0^t  |
R^\varepsilon(u) -R(u)| du +| J^\varepsilon(t) |.
\end{eqnarray*}
Analogously, using that $ E^\varepsilon(u-r) -E(u-r)=0$ for any $u\in (0,r)$, we get:
\begin{eqnarray*}
| E^\varepsilon(t)-E(t) | &\le & \beta \int_0^t | S^\varepsilon(u)-S(u)| du  + \beta \int_0^t |  I^\varepsilon(u)- I(u)| du\\ & &  \qquad+ \int_0^t \frac{1}{K_r} | E^\varepsilon(u-r) -E(u-r)| du
+| J^\varepsilon(t) |\\
&\le & \beta \int_0^t | S^\varepsilon(u)-S(u)| du  + \beta \int_0^t
|  I^\varepsilon(u)- I(u)| du\\ & &  \qquad+ \frac{1}{K_r}  \int_0^t
| E^\varepsilon(u) -E(u)| du +| J^\varepsilon(t) |,
\end{eqnarray*}
\begin{eqnarray*}
|I^\varepsilon(t)-I(t) | &\le &  \frac{1}{K_r} \int_0^t | E^\varepsilon(u) -E(u)| du + \mu  \int_0^t | I^\varepsilon(u)-I(u)| du,\\
|R^\varepsilon(u)-R(u) |&\le &  \mu \ \int_0^t | I^\varepsilon(u)-I(u)| du+
\gamma \int_0^t  | R^\varepsilon(u) -R(u)| du.
\end{eqnarray*}
Putting together these inequalities we obtain the existence of two positive constants $K_1$ and $K_2$ such that
$$
|Z^\varepsilon(t)-Z(t)| \le K_1 | J^\varepsilon(t) |^2 + K_2  \int_0^t  | Z^\varepsilon(u) -Z(u)| du.
$$
Applying clasical Gronwall's lemma, we get for all $t \in [0,T]$:
$$
|Z^\varepsilon(t)-Z(t)| \le K_1 | J^\varepsilon(t) |^2 \exp( K_2  T).
$$
So,
\begin{eqnarray*}
&&P( \Vert Z^\varepsilon-Z \Vert_{\infty,[0,T]} > \rho ) \,\, \le
\,\,P \Big( \Vert J^\varepsilon \Vert_{\infty,[0,T]}^2>
\frac{\rho^2}{K_1 \exp(K_2 T)} \Big)
\\
&= & P \Big( \sup_{t \in [0,T]} |      \int_0^t S^\varepsilon(s) I^\varepsilon(s) dW(s)     |^2> \frac{\rho^2}{\varepsilon^2 K_1 \exp(K_2 T)} \Big)  \\ & \leq &
\exp \Big( - \frac{\rho^2}{\varepsilon^2 K_1 T \exp(K_2 T)} \Big) ,
\end{eqnarray*}
where in the last inequality we have used the exponential martingale inequality and the fact that
$\int_0^T (S^\varepsilon(s) I^\varepsilon(s))^2 ds \le T.$
\hfill$\square$

\subsection{Analysis of the nondelayed system}

Using Lyauponov functionals (see e.g. \cite{HTMW}), we can get a condition for $\varepsilon$ such that the free disease equilibrium is
asymptotically stable for the non delayed stochastic system.

\begin{proposition}

Assume that  $\beta, \mu, \gamma \in (0,1)$; if
\begin{equation*}
\mu > \frac{\beta+\sqrt{\beta^2+2 \varepsilon^2/K_r}}{2}
\end{equation*}
then the free disease equilibrium point is globally asymptotically stable for the non delayed
stochastic system. Note that this condition implies that $\mu>\beta$.
\end{proposition}

\noindent{\bf Proof:}
We prove the stability of the disease-free equilibrium $E_0=(1,0,0,0).$  Using that $S(t)=1-E(t)-I(t)-R(t),$ we can consider that we have a system with three equations.  Putting $u_1=E, u_2=I, u_3=R$ we can consider the linearized system around $(0,0,0)$:
\begin{equation}
\label{SEIRS_nondelay_s}
\left\{
\begin{array}{ll}
du_1(t)= & (\beta u_2(t)  - \frac{1}{K_r} u_1(t) )dt - \varepsilon u_2(t) dW(t)
\\
d u_2(t)= &(\frac{1}{K_r} u_1(t) - \mu u_2(t) ) dt
\\
d u_3(t) = & (\mu u_2(t) - \gamma u_3(t)) dt
\end{array}
\right.
\end{equation}
We denote $u=(u_1,u_2,u_3)$ and we consider the function
$$V(u)= u_1^2 + V_2 u_2^2 + V_3 u_3^2,$$
with $V_2, V_3 >0$. Clearly $V \ge 0$ and $V(0,0,0)=0$. We have
\begin{equation*}
LV  =  2 ( \beta + V_2 \frac{1}{K_r}) u_1 u_2 + 2 V_3 \mu u_2 u_3 -2 \frac{1}{K_r} u_1^2 - (2V_2 \mu -\varepsilon^2 )u_2^ 2 - 2 V_3 \gamma u_3^2 .
\end{equation*}
To get that $LV \le 0$ and using that $2u_1 u_2 \le \lambda_1^2 u_1^ 2 + \frac{1}{\lambda_1^2} u_2^2$ and
$2u_2 u_3 \le \lambda_3^2 u_3^ 2 + \frac{1}{\lambda_3^2} u_2^2$, it is enough to impose that
\begin{eqnarray*}
 -2 \frac{1}{K_r}  + \lambda_1^2 ( \beta + V_2 \frac{1}{K_r})  & \le & 0 \\
 -2V_2 \mu +\varepsilon^2 +\frac{1}{\lambda_1^2}( \beta + V_2 \frac{1}{K_r})+\frac{1}{\lambda_3^2}V_3 \mu& \le & 0 \\
 - 2 V_3 \gamma + \lambda_3^2V_3 \mu & \le & 0
\end{eqnarray*}
Choosing $\lambda_1^2=\frac{2 \frac{1}{K_r} - \alpha_0}{\beta + V_2 \frac{1}{K_r}}>0$ with  $\alpha_0>0$, $\lambda_3^2=\frac{\gamma}{\mu}$ and $V_3$ small enough, choosing $\alpha_0$ as small as we want it suffices to get that
$$
-2V_2 \mu +\varepsilon^2 +\frac{1}{2 \frac{1}{K_r} }( \beta + V_2 \frac{1}{K_r})^2<  0
$$
 i.e.
\begin{equation}\label{negativ}
\frac{1}{K_r^2} V_2^2 + 2 \frac{1}{K_r}  (\beta - 2 \mu) V_2  + \beta^2 + 2 \frac{1}{K_r}  \varepsilon^2 < 0.\end{equation}
Assuming that $2 \mu - \beta >0$, the minimum will be at the point  $V_2 = K_r(2 \mu -\beta)$.
Furthermore, if
\begin{equation}\label{condis}
\mu - \beta - \frac{1}{\mu K_r}\frac{\varepsilon^2}{2} >0.
\end{equation}
which is equivalent, for positive $\mu$, to be
\[
\mu > \frac{\beta+\sqrt{\beta^2+2 \varepsilon^2/K_r}}{2},
\]
then Inequality (\ref{negativ}) holds.
Thus, the proof finishes applying Theorem \ref{lya} with
\begin{eqnarray*}
a(\vert u \vert) &:=& \min \Big( 1, V_2, V_3  \Big) \vert u \vert^ 2,\\
b(\vert u \vert) &:=& \max\Big( 1, V_2, V_3  \Big) \vert u \vert^ 2,\\
c(\vert u \vert) &:=& \min\Big( 2 \frac{1}{K_r}  - \lambda_1^2 ( \beta + V_2 \frac{1}{K_r})  ,
 2V_2 \mu -\varepsilon^2 -\frac{1}{\lambda_1^2}( \beta + V_2 \frac{1}{K_r})-\frac{1}{\lambda_3^2}V_3 \mu,
\\ &  & \qquad
  2 V_3 \gamma - \lambda_3^2V_3 \mu   \Big) \vert u \vert^ 2.
\end{eqnarray*}
\hfill$\square$

\section{Conclusions and future work}\label{con}

In this paper we have analysed the stability of the equilibrium points of a family of SEIRS models.
We consider both deterministic and stochastic models with or without delay proving that
the free disease equilibrium is, under suitable assumptions, always asymptotic stable and that a
similar result for the coexistence equilibrium only holds in some cases.

As a future work we plan to extend these results to models where different time delays
are present in all the equations, to better describe the epidemic models.
For example, if we consider the following general model, for $t>r_{EI}$,
\begin{equation}
\label{SEIRS}
\left\{
\begin{array}{ll}
dS(t) = & -\beta S(t) I(t) dt + \gamma R(t) dt
\\
dE(t) = & \beta S(t) I(t) dt - \beta S(t-r_E) I(t-r_E) dt
\\
dI(t) = & \beta S(t-r_E) I(t-r_E) dt - \beta S(t-r_{EI}) I(t-r_{EI}) dt
\\
dR(t) = & \beta S(t-r_{EI}) I(t-r_{EI}) dt - \gamma R(t) dt
\end{array}
\right.
\end{equation}
where we assume that
any individual remains in the classes E and I, respectively,
for a constant amount of time equal to $r_E$ and $r_I$ and $r_{EI}=r_E+r_I$.


It is easy to see that the possible stability points for such a model are
$(1,0,0,0)$ when $r_I<\frac{1}{\beta}$ and
$(1,0,0,0)$ and
\[
\left( \frac{1}{\beta r_I}, \frac{r_E}{r_I} \frac{\gamma(\beta
r_I-1)}{\beta(\gamma r_{EI}+1)}, \frac{\gamma(\beta
r_I-1)}{\beta(\gamma r_{EI}+1)}, \frac{1}{r_I} \frac{(\beta
r_I-1)}{\beta(\gamma r_{EI}+1)} \right) \
\]
for  $r_I\ge \frac{1}{\beta}$.
In a forthcoming paper we will deal with the study of the stability of the previous equilibrium points,
since the techniques applied in this paper appears not adequate.

\section{Appendix}

In this appendix we recall some well-known results about stability.
First we deal  with the study of the characteristic roots for
deterministic delayed models. Finally we give some results about
stochastic stability using Lyapunov functionals.

A stable steady sate in a deterministic model can become unstable if,
by increasing the delay, a characteristic root changes from having a negative
real part to having positive real part. We will recall here some results about
characteristic functions of order two and three.

\subsection{ Deterministic case: the degree two equation}\label{secdeg2}
Consider the characteristic function of degree two associated to a delayed system
\begin{equation}\label{deg2}
\lambda^2 + a_1 \lambda + a_0 + \big( b_1 \lambda  +b_0 \big)e^{-\lambda r}=0.
\end{equation}
A steady state in this case is stable for $r=0$ if all the roots of
$$
\lambda^2 + (a_1 + b_1) \lambda + ( a_0 + b_0)=0
$$
have negative real part. This occurs if and only if $a_1+b_1>0$ and $a_0+b_0>0$
(by Routh-Hurwitz conditions).  We recall a result about the delayed system (Proposition 2.3 in \cite{For}):

\begin{proposition}\label{Prop2}
A steady state with characteristic equation (\ref{deg2}) is stable in the absence of delay, and becomes unstable with increasing delay if and only if
\begin{enumerate}
\item $a_0 + b_0 >0$ and $a_1+b_1>0$, and

\item $a_0^2 < b_0^2$, or $a_0^2  >b_0^2$, $a_1^2 < b_1^2 + 2a_0$ and $(a_1^2-b_1^2-2a_0)^2>4(a_0^2-b_0^2).$
\end{enumerate}

\end{proposition}

Moreover, as a particular case of the results in Kuang in \cite{Ku} page 74-76, we have that under the hypothesis $a_0+b_0 \not= 0$ and $a_0^2 < b_0^2$, the characteristic function
(\ref{deg2}) has only one imaginary solution $\lambda = i \omega, \omega>0$
\begin{equation}\label{arrel1}
\omega^2= \frac12  \Big( (b_1^2+2a_0-a_1^2) + \big((b_1^2+2a_0-a_1^2)^2-4(a_0^2-b_0^2)\big)^\frac12 \Big).
\end{equation}
Then, the only crossing of the imaginary axis is from left to right as the delay increases. So the stability can only be lost and not regained. Furthermore
the steady state remains asymptotically stable until $r^*=\frac{\theta}{\omega}$, where $ 0 \le \theta <2 \pi$ with
\begin{eqnarray}
\cos \theta & = & - \frac{ a_1 b_1 \omega^2 + (a_0 - \omega^2) b_0 }{ b_1 \omega^2 + b_0^2} \label{eqcos} \\
\sin \theta & = & \frac{ a_1 b_0 \omega - (a_0 - \omega^2) b_1 \omega }{ b_1 \omega^2 + b_0^2}. \label{eqsin}
\end{eqnarray}

\subsection{ Deterministic case: the degree three equation}\label{secdeg3}
We consider now a three degree  general characteristic equation
\begin{equation}\label{deg3}
\lambda^3 + a_2 \lambda^2 + a_1 \lambda + a_0 + \big( b_2 \lambda^2 +b_1 \lambda  +b_0 \big)e^{-\lambda r}=0.
\end{equation}
As in the degree two case, a steady state is stable for $r=0$ if all the roots of
$$
\lambda^3 + (a_2+b_2) \lambda^2 + (a_1 + b_1) \lambda + ( a_0 + b_0)=0
$$
have negative real part. This occurs if and only if $a_2+b_2>0, a_0+b_0>0$ and $(a_2+b_2)(a_1+b_1)-(a_0+b_0)>0$.
Moreover, we will recall a result about the delayed system (Proposition 2.4 in \cite{For}). Set:
\begin{equation}\label{ABC}
A := a_2^2-b_2^2-2a_1, \, B:=a_1^2-b_1^2+2b_2b_0-2a_2a_0 \,\,{ {\rm and }} \, \, C:=a_0^2-b_0^2.
\end{equation}
Then:

\begin{proposition}\label{Prop24}
A steady state with characteristic equation (\ref{deg3}) is stable in the absence of delay, and becomes unstable with increasing delay if and only if
$A, B$ ans $C$ are not all positive and
\begin{enumerate}
\item $a_2+b_2>0, a_0+b_0>0$ and $(a_2+b_2)(a_1+b_1)-(a_0+b_0)>0$, and

\item either $C<0,$  or $C>0,\, A^2-3B>0$ and $4(B^2-3AC)(A^2-3B)-(9C-AB)^2>0.$
\end{enumerate}
\end{proposition}

\subsection{ Stochastic stability}

Consider the $n$ dimensional  stochastic system
$$
dX(t)=f(t,X(t))dt + g(t,X(t))dW_t
$$
where $f(t,x)$ is a function in $\R^n$ defined in $[t_0,+\infty) \times \R^n$, and $g(t,x)$ is a $n \times n$ matrix, $f,g$ are Locally Lipschitz functions in $x$ and $W$ is an $m$-dimensional Wiener process.

Let us denote by $L$ the associated differential operator, defined for a non-negative function $V(t,x) \in C^{1,2} (\R \times \R^n)$ by
$$
LV= \frac{\partial V}{\partial t} + f^T \cdot \frac{\partial V}{\partial x} + \frac12 Tr \big[ g^ T \cdot \frac{\partial^2 V}{\partial x^2} \cdot g \big].$$
Recall that $V$ is called a Lyapunov functional.
The result about stability states as follows:

\begin{theorem}\label{lya}
Suppose that there exist a non-negative function $V(t,x) \in C^{1,2} (\R \times \R^n),$ two continuous function $a,b: \R_+^0 \to \R_+^0$, positive on $R_+$ and a positive constant $K$ such that, for $\vert x \vert <K,$
$$
a(\vert x \vert ) \le V(t,x) \le b(\vert x \vert)$$
holds.
If there exists a continuous function $c: \R_+^0 \to \R_+^0$, positive on $R_+$ such that
$$L V \le  - c(\vert x \vert)$$
holds, then the trivial solution ($X(t)=0$) is globally asymptotically stable.
\end{theorem}

Recall that if $X(t;s,y)$ denotes the solution with initial condition $X(s)=y$ global asymptotic stability means that $\forall \epsilon >0$ and $s \ge t_0$
$$
\lim_{y \to 0} P \Big( \sup_{t \ge t_0} \vert X(t;s,y) \vert \ge \epsilon \Big)=0$$
and
$$
\lim_{y \to 0} P \Big( \lim_{t \to +\infty} \vert X(t;s,y) \vert  = 0 \Big)=0.$$

We refer the reader to \cite{Has} and \cite{Tor} for a complete study of these results.


\end{document}